\documentclass[11pt]{article}
\usepackage{latexsym}
\usepackage{amsfonts}
\usepackage{graphicx}
\usepackage{amsthm}
\usepackage{amsmath}
\usepackage{amssymb}
\usepackage{verbatim}
\usepackage{hyperref}
\usepackage{multirow}
\usepackage{multicol}
\usepackage{enumerate}
\usepackage{tabularx}
\usepackage{lipsum}
\usepackage[arrow,matrix,curve,cmtip,ps]{xy}

\title{Research on Arithmetic}
\author{by Joseph-Louis Lagrange\thanks{Translation by Erik R. Tou, University of Washington Tacoma. Email: \tt{etou@uw.edu}.}}
\date{Published 1775\thanks{\emph{Recherches d'Arithm\'etique}. This article appeared originally in the \emph{Nouveaux M\'emoires de l'Acad\'emie royale des Sciences et Belles-Lettres de Berlin} {\bf 1773} (1775), pp. 265-312. A sequel to this article, \emph{Suite de Recherches d'Arithm\'etique}, appeared in the same journal for the year 1775 (pp. 323-356?).}}

\begin{document}

\newpage
\maketitle

\bigskip\bigskip



This research has for its object the numbers which may be represented by the formula $Bt^2 + Ctu + Du^2$, where $B$, $C$, $D$ are assumed to be given whole numbers, and $t$, $u$ are also whole (though variable) numbers.\footnote{Lagrange used the word ``ind\'etermin\'es'' in French, though here we refer to such ``indeterminate numbers" as variables.} First, I will give the manner by which to find all the different forms whose divisors are the kind of numbers that are susceptible [to this representation]. Next, I will give a method for reducing these forms to the smallest number possible: I will show how one may draw up tables for the practice, and I will note the use of these tables in the research of the divisors of these numbers. Finally, I will prove several Theorems on prime numbers of the same form $Bt^2 + Ctu + Du^2$, of which some are already known, but have not yet been proven, and of which others are entirely new. \\

\begin{center}
\textsc{Note.}
\end{center}

1. We always suppose in the following that all letters designate whole numbers (either positive or negative), and that ordinarily we will represent given numbers by the first letters of the alphabet, and variables by the last [letters of the alphabet]. \\

\begin{center}
\textsc{Observation.}
\end{center}

2. The first-degree formula $Bt+Cu$, where $B$ and $C$ are arbitrary given numbers, and relatively prime, may represent any arbitrary number. But the same is not true for the second-degree formula $Bt^2 + Ctu + Du^2$, since we have proven elsewhere (see the \emph{M\'emoires de l'Acad\'emie} for the years 1767 and 1768) 
that the equation
\[
A \;=\; Bt+Cu
\]
is always resolvable in whole numbers whatever the numbers $A$, $B$, $C$, provided that the last two are relatively prime, but the equation
\[
A \;=\; Bt^2 + Ctu + Du^2
\]
is so only in certain cases, and when certain conditions are placed on the given numbers $A$, $B$, $C$, $D$. We must say the same thing, with even more reason, of third-degree formulas and beyond. \\

\begin{center}
\textsc{Scholium.}
\end{center}

3. There is then a great difference between first-degree formulas and those of higher degrees, the former which can represent all the possible numbers, and the latter which can only represent certain numbers that must be distinguished from all the others by particular characteristics. Very great Geometers have already considered the properties of numbers which can be represented by certain formulas of the second (or higher) degree, such as
$t^2+u^2$, $t^2+2u^2$, $t^2+3u^2$, $t^4+u^4$, $t^8+u^8$, etc. (see the Works of M. de Fermat and the \emph{Novi Commentarii} of St. Petersburg, Vol. 1, 4, 5, 6, 8) 
but I know of no person who has yet treated this matter in a direct and general manner, nor given rules for finding \emph{a priori} the main properties of the numbers which can [be expressed via] given formulas. \\ 

As this subject is a very curious one in Arithmetic, and merits particularly the attention of the Geometers by the great difficulties it contains, I will try to treat it more thoroughly than has yet been done. But at present I will limit myself to second-degree formulas, and I will begin by examining what the form must be of the divisors of those numbers which may be expressed by these kinds of formulas. \\

\begin{center}
\textsc{Theorem I.}
\end{center}

\emph{
4. If the number $A$ is a divisor of a number represented by the formula 
\[
Bt^2 + Ctu + Du^2,
\]
supposing $t$ and $u$ are relatively prime, I say that this number $A$ will necessarily have the form
\[
A \;=\; Ls^2+Msx+Nx^2,
\] 
where we will have 
\[
4LN - M^2 \;=\; 4BD - C^2,
\]
$s$ and $x$ also being relatively prime.
} \\

Let $a$ be the quotient obtained from the division of $Bt^2 + Ctu + Du^2$ by $A$, so that we have 
\[
Aa\;=\; Bt^2 + Ctu + Du^2,
\]
and let $b$ be the greatest common divisor of $a$ and $u$ (if $a$ and $u$ are relatively prime, we will have $b=1$), so that by writing $a = bc$ and $u = bs$, $c$ and $s$ are relatively prime. Therefore we have 
\[
Abc \;=\; Bt^2+Cbts+Db^2s^2;
\]
consequently $Bt^2$ will be divisible by $b$. But, $t$ and $u$ being relatively prime (by hypothesis), $t$ will also be relatively prime to $b$, which is a divisor of $u$. Therefore it must be that $B$ is divisible by $b$. So we will have $B = Eb$, and, the equation being divisible by $b$, it will become 
\[
Ac \;=\; Et^2 + Cts + Dbs^2.
\]
Now, since $c$ and $s$ are relatively prime, we may assume (by the preceding Observation) that $t = \theta s + cx$, which, being substituted, will give
\[
Ac \;=\; (E\theta^2+C\theta+Db)s^2 + (2E\theta c+Cc)sx + Ec^2x^2,
\]
so that it will be necessary that the number $(E\theta^2+C\theta+Db)s^2$ is divisible by $c$. And since $c$ and $s$ are relatively prime, it must be that $E\theta^2+C\theta+Db$ is divisible by $c$; therefore dividing the whole equation by $c$, and setting 
\begin{eqnarray*}
L &=& \frac{E\theta^2+C\theta+Db}{c}, \\
M &=& 2E\theta+C, \\
N &=& Ec,
\end{eqnarray*}
we will have $A = Ls^2+Msx+Nx^2$. \\ 

Now $4LN - M^2$ will be equal to $4E(E\theta^2+C\theta+Db) - (2E\theta+C)^2 = 4EDb-C^2 = 4BD-C^2$, because $B = Eb$. Therefore, etc. \\

Now, since $t$ and $u$ are relatively prime (by hypothesis), $t$ and $s$ will be as well, because $u = bs$. But if $x$ and $s$ are not relatively prime it is clear that $t$ should be divisible by their greatest common divisor, because $t = \theta s+cx$. Since this cannot be, it follows that $x$ and $s$ will necessarily be relatively prime whenever $t$ and $u$ are [relatively prime]. \\

\begin{center}
\textsc{Theorem II.}
\end{center}

5. Every second-degree formula such as this one
\[
Ls^2+Msx+Nx^2,
\]
in which $M$ is larger than $L$ or $N$ (ignoring the signs of these quantities), may be transformed into another of the same degree, as
\[
L'{s'}^2+M's'x'+N'{x'}^2,
\]
for which we will have 
\[
4L'N' - {M'}^2 \;=\; 4LN - M^2,
\]
and where $M'$ will be less than $M$. \\

Suppose, for example, that $M > L$. We will have $s = mx+s'$, and the proposed formula will become
\[
(Lm^2+Mm+N)x^2 + (2Lm + M)xs' + L{s'}^2,
\]
or, by changing $x$ to $x'$,
\[
L'{s'}^2+M's'x' + N'{x'}^2,
\]
where we will have 
\begin{eqnarray*}
L' &=& L, \\
M' &=& 2Lm+M, \\
N' &=& Lm^2+Mm+N.
\end{eqnarray*}
So first we will have, regardless of $m$,
\[
4L'N'-{M'}^2 = 4L(Lm^2+Mm+N)-(2Lm+m)^2 = 4LN-M^2.
\]
Now, since $L$ is less than $M$ (by hypothesis), it is clear that we may determine the number $m$ so that [$M'=$ ] $2Lm+M$ becomes less than $M$; therefore, etc. \\

\begin{center}
\textsc{Corollary 1.}
\end{center}

6. Therefore, if in the transformed $L'{s'}^2+M's'x'+N'{x'}^2$ one of the numbers $L$ or $N'$ is less than $M'$, we will be able to obtain another transformation such that
\[
L''{s''}^2 + M''s''x''+N''{x''}^2,
\] 
in which we will likewise have $4L''N''-{M''}^2 = 4L'N' - {M'}^2 = 4LN-M^2$ and where $M''$ will be less than $M'$; and so on. Therefore, since the series\footnote{Lagrange uses \emph{s\'erie} here, though a more modern rendering would refer to this as a sequence.} of numbers $M$, $M'$, $M''$, etc. cannot go to infinity (because the numbers must all be whole, and decreasing from one to the next) we must necessarily arrive at a transformation which I will represent as
\[
Py^2 + Qyz + Rz^2,
\]
in which $Q$ will not be greater than $P$, nor $R$, and where we will have $4PR-Q^2 = 4LN - M^2$.

\begin{center}
\textsc{Corollary 2.}
\end{center}

7. If the numbers $s$ and $x$ from the proposed formula are relatively prime, it is clear that the numbers $s'$ and $x'$ from the transformation will also be relatively prime; because if they were not it would follow necessarily (because $x' = x$ and $s = mx+s'$) that $s$ was divisible by the greatest common divisor of $s'$ and $x$. \\

Therefore the numbers $s''$ and $x''$ from the second transformation will also be (for the same reason) relatively prime, and so on, from which on may conclude that the numbers $y$ and $z$ from the last transformation will necessarily be relatively prime, if the numbers $s$ and $x$ are [relatively prime]. \\

\begin{center}
\textsc{Theorem III.}
\end{center}

8. If $A$ is a divisor of a number of the form 
\[
Bt^2+Ctu+Du^2,
\]
with $t$ and $u$ being relatively prime, I say that the number $A$ will necessarily have the form 
\[
Py^2+Qyz+Rz^2,
\]
with $y$ and $z$ also being relatively prime, and $P$, $Q$, $R$ being such that 
\[
4PR-Q^2 \;=\; 4BD - C^2,
\]
$Q$ not being greater than $P$ or $R$ (ignoring the signs of $P$, $Q$, and $R$). \\

The proof of this Theorem follows naturally from the two preceding theorems, and their Corollaries. \\

\begin{center}
\textsc{Corollary 1.}
\end{center}

9. If $4BD-C^2$ is positive, it must be that $4PR$ is also positive; so, because $P\ge Q$ and $R \ge Q$, it is clear that $4PR$ will also be $\ge 4Q^2$, and consequently $4PR-Q^2 \ge 3Q^2$. Therefore, we will also have $4BD-C^2 \ge 3Q^2$, and from this [it follows that]
\[
Q \;\le\; \sqrt{\frac{4BD-C^2}{3}}.
\]

\begin{center}
\textsc{Corollary 2.}
\end{center}

10. Now suppose $4BD-C^2$ is negative, so that $C^2-4BD$ is positive; since $Q$ is not greater than $P$ or $R$ we  will therefore have that the case $Q^2-4PR > 0$ cannot happen unless $4PR$ is a negative number; so $-4PR$ will be a positive number $\ge 4Q^2$, because $P\ge Q$ and $R\ge Q$. So then $Q^2-4PR$ will be $\ge 5Q^2$, and consequently $C^2-4BD$ will also be $\ge 5Q^2$; therefore it must be that 
\[
Q \;\le\; \sqrt{\frac{C^2-4BD}{5}}.
\]

\begin{center}
\textsc{Corollary 3.}
\end{center}

11. Therefore, since $Q$ is a whole number, we can only take for $Q$ positive or negative whole numbers which do not surpass these discovered limits, also understanding zero to be among these whole numbers. From this we see that $Q$ can only have a fixed number of different values. \\ 

Furthermore, it is clear that for the equation $4PR-Q^2 = 4BD-C^2$ to subsist in whole numbers, it must be that $Q$ is even or odd, consequently that $C$ will be even or odd, which even further limits the number of values of $Q$. \\

Knowing $Q$, we will easily find $P$ and $R$ by the same equation; so because $PR = \frac{4BD-C^2+Q^2}{4}$ it is clear that we can only take for $P$ and $R$ the factors of the whole number $\frac{Q^2+4BD-C^2}{4}$, taking care to reject those for which one or both would be greater than $Q$. \\

\begin{center}
\textsc{Problem I.}
\end{center}

12. \emph{Find all the possible forms for the divisors of those numbers which are represented by the second-degree formula}
\[
Bt^2+Ctu+Du^2,
\]
\emph{with $t$ and $u$ being relatively prime.} \\

This is obvious since we have just demonstrated above that each divisor of the proposed formula is reducible to the form
\[
Py^2+Qyz+Rz^2,
\]
with $y$ and $z$ being relatively prime. So, the difficulty is reduced to finding the values of the coefficients $P$, $Q$, $R$ when $B$, $C$, and $D$ are given. \\

To this end, I distinguish between two cases: one in which the number $4BD-C^2$ is positive, and the other when this number is negative. \\

First Case. Let $4BD-C^2 = K$, with $K$ designating a positive number. First, we will determine $Q$ from these conditions: from whether $Q$ is even or odd it follows what $K$ will be, and if it does not exceed the number $\pm\sqrt{\frac{K}{3}}$. 
Then we will determine $P$ and $R$ by these conditions: that $P$ and $R$ are two factors of the number $\frac{K+Q^2}{4}$, and that each of these factors is not less than $Q$ (Articles 9 and 11). \\

Second Case. Let $4BD-C^2 = -K$. We will determine $Q$ by these conditions: whether $Q$ is even or odd it follows what $K$ will be, and if it does not exceed the number $\pm\sqrt{\frac{K}{5}}$. 
Next we will determine the corresponding values of $P$ and $R$ by these conditions: that $P$ and $R$ are two factors of the number $\frac{Q^2-K}{4}$, and that each of them is not less than $Q$ (Articles 10 and 11). \\

\begin{center}
\textsc{Remark 1.}
\end{center}

13. If we have $4BD-C^2 = 0$, then with $K$ being $=0$ we could only take $Q = 0$, and then we will also have $PR = 0$, so that one of the numbers $P$ or $R$ will be zero (and the other can be anything we would like). But it is necessary to remark that in this case the formula $Bt^2+Ctu+Du^2$ reduces to $\frac{(2Bt+Cu)^2}{4}$, so since $2Bt+Cu$ may represent any arbitrary number (Article 2), the divisors of the proposed formula may also be arbitrary. \\

\begin{center}
\textsc{Remark 2.}
\end{center}

14. The same thing must take place, in general, when the formula $Bt^2+Ctu+Du^2$ is the product of two rational, first-degree formulas such as $at+bu$ and $ct+du$, each of which may represent arbitrary numbers (Article 2). This is what occurs when $4BD-C^2$ is equal to a negative square. So, supposing $4BD-C^2 = -H^2$, we have 
\[
Bt^2+Ctu+Du^2 \;=\; \frac{(2Bt+(C+H)u)(2Bt+(C-H)u)}{4B}.
\]
Now, though in this case any number can be a divisor of the formula in question, if we search the formulas of the divisors in the preceding Problem, we will find them to be as in the other case, so that we must conclude that these formulas contain all the possible numbers. \\

For the rest we have $4PR-Q^2 = 4BD-C^2 = -H^2$, and it is clear that the general formula of the divisors $Py^2+Qyz+Rz^2$ will also be resolvable in two rational, first-degree formulas. \\

\begin{center}
\textsc{Remark 3.}
\end{center}

15. It is remarkable that the formulas for the divisors do not depend on the value of $K$, which is to say, the number $4BD-C^2$. But it is easy to see the reason by noting that the formula $Bt^2+Ctu+Du^2$ may reduce to $\frac{(2Bt+Cu)^2+(4BD-C^2)u^2}{4B}$, so that the divisors of the formula $Bt^2+Ctu+Du^2$ may also be regarded as divisors of the more simple formula $x^2\pm Ku^2$. \\

From this it follows that it suffices to consider the formulas of this last type, and for that we will add the following Problem, which may be regarded as a special case of the preceding one, but which fundamentally has the same generality. \\

\begin{center}
\textsc{Problem II.}
\end{center}

16. \emph{To find the possible formulas for the divisors of numbers of the form $t^2\pm au^2$, $a$ being a given arbitrary positive integer, and $t$ and $u$ being variable, relatively prime numbers.}
\\

First, consider the formula $t^2+au^2$. Comparing to the general formula from Problem I, we have $B=1$, $C=0$, $D=a$, thus $K = 4a$; therefore $Q$ will be even and it will not be greater than $\pm\sqrt{4a}{3}$; also making $Q = \pm 2q$ and regarding $q$ as positive, it must be that $q$ is not $> \sqrt{\frac{a}{3}}$; then we have $PR = \frac{4a+4q^2}{4} = a+q^2$. So then if $p$ and $r$ denote two factors of $a+q^2$, neither of which is less than $2q$, we will have 
\[
py^2\pm 2qyz + rz^2
\]
for the general formula of the divisors of $t^2+au^2$. \\

It is worth noting that since $pr=a+q^2$, $p$ and $r$ must have the same sign, and it is clear that they must be positive so that $py^2\pm 2qyz+rz^2$ can represent positive numbers. \\

Furthermore, since this formula does not change form (in terms of $y$) when putting $p$ in place of $r$, it will not be necessary to take take successively for $p$ each of the factors of $a+q^2$, and for $r$ all the corresponding factors. That is why, in each pair of factors of $a+q^2$, it will suffice to always take $p$ to be the smallest and $r$ to be the largest. This is how we will use it 
in what follows. \\

Second, consider the formula $t^2-au^2$. We have $B=1$, $C=0$, $D=-a$, thus $K=4a$, as above; this is why $Q = \pm 2q$ will be the same, and it must be that $q$ is not $> \sqrt{\frac{a}{5}}$; so we will have $PR = q^2-a$. So then if we designate by $p$ and $r$ two factors of $a-q^2$, neither of which is less than $2q$, we will have $P=p$, $R=-r$, or $P=-p$, $R=r$. These will give the two formulas
\[
\begin{cases} & py^2\pm 2qyz-rz^2 \\ -\!\!\!\!\!\!\!\!\! & py^2\pm 2qyz+ry^2 \end{cases} 
\]
for the divisors of $t^2-au^2$. We will find the same thing for the formula $au^2-t^2$. \\

As for the numbers $p$ and $r$, we take both to be positive, and suppose always that $p$ is the smaller of the two factors of $a-q^2$, and $r$ is the larger (as we have said until now). Now, it is clear that by changing the signs of $p$ and $r$, or putting one of these numbers in place of the other, we will not produce any new formulas. \\

\begin{center}
\textsc{Corollary.}
\end{center}

17. If we multiply the formula $py^2\pm 2qyz + rz^2$ by $p$, it can be put into the form $(py\pm qz)^2+(pr-q^2)z^2$, which is to say (because $pr = a+q^2$) the form $(py\pm qz)^2 + az^2$, which is the same as the formula $t^2+au^2$. Then it follows that every divisor of a number of the form $t^2+au^2$ will also necessarily be of the same form if $p$ does not have values other than unity, or will become unity after being multiplied by one of the values of $p$ (if there are several). 
We will prove the same, that the formulas $py^2\pm 2qyz-rz^2$ and $-py^2\pm 2qyz+rz^2$ being multiplied by $p$ will become (because $pr = a-q^2$) $(py\pm qz)^2-az^2$ and $-(py\pm qz)^2+az^2$. So then every divisor of a number of the form $t^2-au^2$ or $au^2-t^2$ will be necessarily of the one or the other of these two forms if $p$ is not unity, or else it will always become unity after being multiplied by one of the values of $p$ (if there is more than one). \\

18. \emph{Theorems on the divisors of the numbers of the form
\[
t^2+au^2,
\]
$t$ and $u$ being relatively prime.
}

\begin{center}
\textsc{I.}
\end{center}

Let $a=1$, thus $q$ is not $>\sqrt{\frac{1}{3}}$; so $q=0$ and $pr=1$; thus $p=1$ and $r=1$. \\

Therefore the divisors of the numbers of the form $t^2+u^2$ are necessarily contained in the formula $y^2+z^2$, which is to say that every divisor of a number equal to the sum of two squares is also the sum of two squares. \\

\begin{center}
\textsc{II.}
\end{center}

Let $a=2$, thus $q$ is not $>\sqrt{\frac{2}{3}}$; so $q=0$ and $pr=2$; thus $p=1$ and $r=2$. \\

Therefore the divisors of the numbers of the form $t^2+2u^2$ are contained in the formula $y^2+2z^2$, which is to say that every divisor of a number equal to the sum of a square and a doubled square is also the sum of a square and a doubled square. \\

\begin{center}
\textsc{III.}
\end{center}

Let $a=3$, thus $q$ is not $>\sqrt{\frac{3}{3}} =1$; so $q=0$ or $=1$. Making $q=0$, we will have $pr=3$, so $p=1$ and $r=3$. Next making $q=1$, we have $pr = 3+1=4$. Thus, since neither $p$ nor $r$ can be $<2q$, we must have $p=2$ and $r=2$. \\

Therefore the divisors of the numbers of the form $t^2+3u^2$ will be contained in the two formulas $y^2+3z^2$ and $2y^2\pm 2yz + 2z^2$. Now, since the second of these formulas can only belong to the even numbers, all terms being divisible by $2$, it follows that every odd divisor of $t^2+3u^2$ will be contained necessarily in the formula $y^2+3z^2$. In other words, every odd divisor of a number which is the sum of a square and a tripled square (which are relatively prime) is also the sum of a square and a tripled square. \\ 

Besides, 
since it suffices to consider the odd divisors, in what follows we will always ignore those formulas which would admit only even divisors.

\begin{center}
\textsc{IV.}
\end{center}

Let $a=4$, thus $q$ is not $> \sqrt{\frac{4}{3}}$; so $q = 0$ or $=1$. Making $q=0$, we have $pr=4$, so $p=1$ and $r=4$ (and we reject the values $p=2$ and $r=2$ because they are both even). Making $q=1$, we have $pr=5$; thus $p=1$ and $r=5$, which must be rejected since $p$ will be $<2q$. \\

Therefore the odd divisors of the numbers of the form $t^2+4u^2$ will also be of the form $y^2+4z^2$. \\

\begin{center}
\textsc{V.}
\end{center}

Let $a=5$, thus $q$ is not $>\sqrt{\frac{5}{3}}$; so $q=0$ or $=1$. Making $q=0$, we have $pr=5$, so $p=1$ and $r=5$. Making $q=1$, we have $pr=6$; thus $p=2$ and $r=3$. \\

Therefore the divisors of the numbers of the form $t^2+5u^2$ are necessarily of one or the other of these two forms: $y^2+5z^2$ or $2y^2\pm 2yz+3z^2$. So these divisors themselves or their doubles are always of the form $t^2+5u^2$ (see Section 17). 
\\

\begin{center}
\textsc{VI.}
\end{center}

Let $a=6$, thus $q$ is not $>\sqrt{\frac{6}{3}}$; so $q=0$ or $=1$. Making $q=0$, we will have $pr=6$, so either $p=1$ and $r=6$ or $p=2$ and $r=3$. Making $q=1$, we will have $pr=7$; thus $p=1$ and $r=7$, which must be rejected since $p$ will be $<2q$. \\

Therefore the divisors of the numbers of the form $t^2+6u^2$ will have one or the other of the forms $y^2+6z^2$ and $2y^2+3z^2$, so that these divisors themselves or their doubles will be of the form $t^2+6u^2$. 

\begin{center}
\textsc{VII.}
\end{center}

Let $a=7$, thus $q$ is not $>\sqrt{\frac{7}{3}}$; so $q=0$ or $=1$. Making $q=0$ we will have $pr=7$, so $p=1$ and $r=7$. Making $q=1$, we will have $pr=8$; thus $p=2$ and $r=4$, which can only admit 
even divisors. \\

Therefore the odd divisors of the numbers of the form $t^2+7u^2$ will also be necessarily of the form $y^2+7z^2$. \\   

\begin{center}
\textsc{VIII.}
\end{center}

Let $a=8$, thus $q$ is not $>\sqrt{\frac{8}{3}}$; so $q=0$ or $=1$. Making $q=0$ we will have $pr=8$, so $p=1$ and $r=8$. We will reject the values $p=2$, $r=4$ as they may belong only to even divisors. Next, making $q=1$ we will have $pr=9$; thus $p=3$ and $r=3$. \\

Therefore the divisors of the numbers of the form $t^2+8u^2$ are of the one or the other of the forms $y^2+8z^2$ or $3y^2\pm 2yz + 3z^2$, so that these divisors themselves, or their triples, will be always of the same form $t^2+8u^2$. \\

\begin{center}
\textsc{IX.}
\end{center}

Let $a=9$, thus $q$ is not $>\sqrt{\frac{9}{3}}$; so $q=0$ or $=1$. Making $q=0$, we will have $pr=9$; thus either $p=1$ and $r=9$ or $p=3$ and $r=3$. Making $q=1$, we will have $pr=10$; thus $p=2$ and $r=5$. \\

Therefore the divisors of the numbers of the form $t^2+9u^2$ are necessarily of one of the three forms $y^2+9z^2$, $3y^2+3z^2$, $2y^2\pm 2yz+5z^2$; so that these divisors themselves, or their doubles or their triples, will always be able to relate to the same form $t^2+9u^2$. 
\\

\begin{center}
\textsc{X.}
\end{center}

Let $a=10$, thus $q$ is not $>\sqrt{\frac{10}{3}}$; so $q=0$ or $=1$. Making $q=0$ we will have $pr=10$; thus either $p=1$ and $r=10$, or $p=2$ and $r=5$. Making $q=1$, we will have $pr=11$; thus $p=1$ and $r=11$, which is not possible because then $p$ would be $<2q$. \\

Therefore the divisors of the numbers of the form $t^2+10u^2$ are always of one of the forms $y^2+10z^2$ or $2y^2+5z^2$, so that these divisors themselves or their doubles will be necessarily of the same form $t^2+10u^2$. \\ 

\begin{center}
\textsc{XI.}
\end{center}

Let $a=11$, thus $q$ is not $>\sqrt{\frac{11}{3}}$; so $q=0$ or $=1$. Making $q=0$, we will have $pr=11$; thus $p=1$ and $r=11$. Making $q=1$, we will have $pr=12$; thus $p=3$ and $r=4$ (the values $p=2$ and $r=6$ are rejected because they would only admit even divisors). \\

Therefore the divisors of the numbers of the form $t^2=11u^2$ are of one or the other of the forms $y^2+11z^2$ or $3y^2\pm 2yz+4z^2$; so that these divisors themselves or their triples will be always of the same form $t^2+11u^2$. \\

\begin{center}
\textsc{XII.}
\end{center}

Let $a=12$, thus $q$ is not $>\sqrt{\frac{12}{3}} = 2$; so $q=0$, $=1$, or $=2$. Making $q=0$, we will have $pr=12$; thus either $p=1$ and $r=12$ 
or $p=3$ and $r=4$ (rejecting the values $p=2$ and $r=6$, 
which would only admit even divisors). Making $q=1$, we will have $pr=13$; thus $p=1$ and $r=13$, which must be rejected because $p$ will be $< 2q$. Making $q=2$, we will have $pr=12+4 = 16$; thus $p=4$ and $r=4$ (because $q=2$, $p$ cannot be $<4$) which must be rejected if we only consider odd divisors. \\

Therefore the odd divisors of the numbers of the form $t^2+12u^2$ will be of one or the other of these forms $y^2+12z^2$ or $3y^2+4z^2$; so that these divisors themselves or their triples will be of the same form $t^2+12u^2$. \\

We will not extend these this research any further, especially since the examples which we just gave are more than sufficient to show the application of our methods and to put on the path those who will want to use it to discover new theorems on the form of the divisors of the numbers $t^2+au^2$. \\

\begin{center}
\textsc{Remark.}
\end{center}

19. These three initial Theorems have long been known by the Geometers, and are due (I believe) to Mr. Fermat, though Mr. Euler is the first to prove them. We can see the proofs of them in Volumes IV, VI, and VIII of the \emph{Novi Commentarii} of [the St.] Petersburg [Academy].\footnote{This appears to refer to the following articles by Euler: \emph{De numeris, qui sunt aggregata duorum quadratorum
} [E228], \emph{Solutio generalis quorundam problematum Diophanteorum, quae vulgo nonnisi solutiones speciales admittere videntur} and \emph{Specimen de usu observationum in mathesi pura} [E255-256], \emph{Theoremata arithmetica nova methodo demonstrata} and \emph{Supplementum quorundam theorematum arithmeticorum, quae in nonnullis demonstrationibus supponuntur} [E271-272].} His method is totally different from mine, and anyway it is applicable only to the cases where the number $a$ does not surpass $3$; this is what may have prevented the great Geometer from further pursuing his research on the subject. \\

20. \emph{Theorems on the divisors of the numbers
\[
t^2-au^2 \;\;\; \text{or} \;\;\; au^2-t^2,
\]
$t$ and $u$ being relatively prime.
}

\begin{center}
\textsc{I.}
\end{center}

Let $a=1$, thus $q$ is not $>\sqrt{\frac{1}{5}}$; so $q=0$ and $pr=1$; thus $p=1$ and $r=1$. \\

Therefore the divisors of the numbers of the form $t^2-u^2$ will be of the form $y^2-z^2$; consequently (Section 14) every number is reducible to the form $y^2-z^2$; this is what was done elsewhere. \\

\begin{center}
\textsc{II.}
\end{center}

Let $a=2$, thus $q$ is not $>\sqrt{\frac{2}{5}}$; so $q=0$ and $pr=2$; thus $p=1$ and $r=2$, so that the forms of the divisors of $t^2-2u^2$ or $2u^2-t^2$ will be $y^2-2z^2$ or $2z^2-y^2$. But note that these two forms amount to the same thing: making $y = y'+2z'$ and $z=y'+z'$ (which gives $y'=2z-y$ and $z'=y-z$, so consequently $y'$ and $z'$ are whole numbers) the formula $y^2-2z^2$ becomes $2z^2-y^2$. \\

Therefore the divisors of the numbers of the form $t^2-2u^2$ or $2u^2-t^2$ are necessarily of the one or the other of the forms $y^2-2z^2$ and $2z^2-y^2$. \\

\begin{center}
\textsc{III.}
\end{center}

Let $a=3$, thus $q$ is not $>\sqrt{\frac{3}{5}}$; so $q=0$ and $pr=3$; thus $p=1$ and $r=3$. \\

Therefore the divisors of the numbers of the form $t^2-3u^2$ or $3u^2-t^2$ are of the one or the other of the forms $y^2-3z^2$ and $3z^2-y^2$. 
\\

\begin{center}
\textsc{IV.}
\end{center}

Let $a=4$, thus $q$ is not $>\sqrt{\frac{4}{5}}$; so $q=0$ and $pr=4$; thus either $p=1$ and $r=4$, or $p=2$ and $r=2$. \\

Therefore, the divisors of the numbers of the form $t^2-4u^2$ or $4u^2-t^2$ will be necessarily contained in the formulas $y^2-4z^2$, $4z^2-y^2$, $2y^2-2z^2$; consequently (Section 14) any arbitrary number will have one of these forms. \\

For the rest we may ignore the forms which would only admit even divisors, such as $2y^2-2z^2$; so we reject in the following, as we have done up until now, the values of $p$ and $r$ which are simultaneously even. \\

\begin{center}
\textsc{V.}
\end{center}

Let $a=5$, thus $q$ is not $>\sqrt{\frac{5}{5}}=1$; so $q=0$ or $=1$. Making $q=0$, we have $pr=5$, thus $p=1$ and $r=5$. Making $q=1$, we would have $pr=4$, so because $p$ and $r$ cannot be $<2q$ we can only do $p=2$ and $r=2$ (but we reject these values because they are always both even). So we will not have that these two forms of divisors $y^2-5z^2$ and $5z^2-y^2$, 
which anyhow reduce to the same thing, since we may agree that in making $y=2y'+5z'$ and $z=y'+2z'$ (which gives $z'=y-2z$ and $y'=5z-2y$, and consequently $y'$ and $z'$ are whole numbers) in the formula $y^2-5z^2$ which, by these substitutions, will become $5{z'}^2-{y'}^2$. \\

Therefore, the odd divisors of the numbers of the form $t^2-5u^2$ or $5u^2-t^2$ are in the two forms $y^2-5z^2$ and $5z^2-y^2$, respectively. \\

\begin{center}
\textsc{VI.}
\end{center}

Let $a=6$, thus $q$ is not $>\sqrt{\frac{6}{5}}$; so $q=0$ or $=1$. Making $q=0$, we will have $pr=6$, so either $p=1$ and $r=6$, or $p=2$ and $r=3$. Making next $q=1$, we will have $pr=5$, which would only give $p=1$ and $r=5$, values which are not admissible since $p$ would be $<2q$. So the formulas of the divisors of the numbers of the form $t^2-6u^2$ or $6u^2-t^2$ will be $y^2-6z^2$, $6z^2-y^2$, $2y^2-3z^2$, $3y^2-2z^2$. But I observe that these last two reduce to the first two by making $2y+3z=y'$, $y+z=z'$, which gives $y=3z'-y'$, $z=y'-2z'$, and consequently $2y^2-3z^2 = 6{z'}^2-{y'}^2$, $3z^2-2y^2 = {y'}^2-6{z'}^2$. \\

Therefore, the divisors of the numbers of the form $t^2-6u^2$ or $6u^2-t^2$ will always also have one or the other of these forms. \\

\begin{center}
\textsc{VII.}
\end{center}

Let $a=7$, thus $q$ is not $>\sqrt{\frac{7}{5}}$; so $q=0$ or $=1$. Making $q=0$, we will have $pr=7$, so $p=1$ and $r=7$. Making $q=1$, we will have $pr=6$, thus $p=2$ and $r=3$. So the formulas of the divisors of $t^2-7u^2$ will be $y^2-7z^2$ $2y^2\pm 2yz-7z^2$, and their opposites $7z^2-y^2$, $7z^2\pm 2yz-2y^2$. But I note here that the first two of these formulas become the same thing, as well as the last two; now making $y=y'-2z'$ and $\pm z = y'-3z'$ (which gives $y'=3y\mp 2z$ and $z'=y\mp z$, which is to say that $y'$ and $z'$ are whole numbers) the formula $2y^2\pm 2yz-3z^2$ will become ${y'}^2-y{z'}^2$, and the formula $3z^2\mp 2yz-2y^2$ will become the same $7{z'}^2-{y'}^2$.  \\

From this it follows that the divisors of the numbers of the form $t^2-7u^2$ or $7u^2-t^2$ also will be necessarily  of the forms $y^2-7z^2$ or $7z^2-y^2$. \\

\begin{center}
\textsc{VIII.}
\end{center}

Let $a=8$, thus $q$ is not $>\sqrt{\frac{8}{5}}$; so $q=0$ or $=1$. Making $q=0$, we will have $pr=8$, so either $p=1$ and $r=8$, or $p=2$ and $r=4$. But these last values may be rejected because they are always both even. Making then $q=1$, we will have $pr=7$, thus $p=1$ and $r=7$, which would only give $p=1$ and $r=7$, 
values which are not admissible since $p$ would be $<2q$. \\

Therefore, the odd divisors of the numbers of the form $t^2-8u^2$ or $8u^2-t^2$ will be of the one or the other of the two forms $y^2-8z^2$ or $8z^2-y^2$. \\

\begin{center}
\textsc{IX.}
\end{center}

Let $a=9$, thus $q$ is not $>\sqrt{\frac{9}{5}}$; so $q=0$ or $=1$. Making $q=0$, we will have $pr=9$, so either $p=1$ and $r=9$, or $p=3$ and $r=3$. Making $q=1$, we will have $pr=8$, which, because $p$ is not $<2q$, will give $p=2$ and $r=4$, values which we may reject because they are both even. \\

Therefore the odd divisors of the numbers of the form $t^2-9u^2$ or $9u^2-t^2$ will always be of one of the forms $y^2-9z^2$, $9z^2-y^2$, $3y^2-3z^2$; consequently (Section 14) any arbitrary odd number will be reducible to one of these forms. \\

\begin{center}
\textsc{X.}
\end{center}

Let $a=10$, thus $q$ is not $>\sqrt{\frac{10}{5}} = \sqrt{2}$; so $q=0$ or $=1$. Making $q=0$, we will have $pr=10$, so either $p=1$ and $r=10$, or $p=2$ and $r=5$. Making $q=1$, we will have $pr=9$, thus $p=3$ and $r=3$; so that the forms of the divisors of $t^2-10t^2$ 
will be $y^2-10z^2$, $10z^2-y^2$, $2y^2-5z^2$, $5z^2-2y^2$, and $3y^2\pm 2yz-3z^2$. Now I remark first that this last formula may be reduced to these two here, $2{y'}^2-5{z'}^2$ and $5{z'}^2-2{y'}^2$, making $\pm y = y'+z'$ and $z=y'+2z'$, or $\pm y=y'+2z'$ and $z=y'+z'$, which always gives whole numbers for $y'$ and $z'$. I remark next that the two forms $y^2-10z^2$ and $10z^2-y^2$ may also be reduced to the same by making in the first $y=10z'+3y'$ and $z=3z'+y'$, which will transform it into $10{z'}^2-{y'}^2$; and as to the numbers $y'$ and $z'$ it is clear that they will always be whole, since we will have $z'=y-3z$ and $y'=10z-3y$. \\

From there, I conclude that the divisors of the numbers of the form $t^2-10u^2$ or $10u^2-t^2$ will always be of the one or the other of these two forms $y^2-10z^2$ or $2y^2-5z^2$, as well as these, $10z^2-y^2$ or $5z^2-2y^2$. \\

\begin{center}
\textsc{XI.}
\end{center}

Let $a=11$, thus $q$ is not $>\sqrt{\frac{11}{5}}$; so $q=0$ or $q=1$. Making $q=0$, we will have $pr=11$, so $p=1$ and $r=11$. Making $q=1$, we have $pr=10$, so $p=2$ and $r=5$. So in this case, the forms of the divisors will be $y^2-11z^2$, $11y^2-z^2$, $2y^2\pm 2yz-5z^2$, $5z^2\pm 2yz -2y^2$. But I note that these last two formulas may be reduced to the first two; so in making $\pm y=y'+4z'$, $z=y'+3z'$ (which gives $z'=\pm y-z$ and $y'=4z\mp 3y$, and consequently $y'$ and $z'$ are always whole numbers) the formula $2y^2\pm 2yz -5z^2$ becomes $11{z'}^2-{y'}^2$, and the formula $5z^2\mp 2yz-2y^2$ becomes the same ${y'}^2-11{z'}^2$. \\

From this it follows that the divisors of the numbers of the form $t^2-11u^2$ or $11u^2-t^2$ will always be of the one or the other of the forms $y^2-11z^2$ or $11z^2-y^2$. \\

\begin{center}
\textsc{XII.}
\end{center}

Let $a=12$, so $q$ is not $>\sqrt{\frac{12}{5}}$; so $q=0$ or $=1$. Making $q=0$, we will have $pr=12$, so either $p=1$ and $r=12$, or $p=3$ and $r=4$, rejecting the even values $p=2$ and $r=6$. Next making $q=1$, we would have $pr=11$, thus $p=1$ and $r=11$, values which are not admissible because $p$ would be $<2q$. So we will only have the formulas $y^2-12z^2$, 
$12z^2-y^2$, $3y^2-4z^2$, $4z^2-3y^2$, for which I remark that these last two are reducible to the first two, by making $y=4y'+z'$ and $z=3y'+z'$ (which gives $y'=y-z$ and $z'=4z-3y$, and consequently $y'$ and $z'$ are whole numbers). \\

From this we may conclude that the odd divisors of the numbers of the form $t^2-12u^2$ or $12u^2-t^2$ will always be of the one or the other of the two forms $y^2-12z^2$ or $12z^2-y^2$, as well as these two, $3y^2-4z^2$ or $4z^2-3y^2$. \\

\begin{center}
\textsc{Remark.}
\end{center}

21. Such is the method which it is necessary to follow in order to find the formulas for divisors of numbers of the form $t^2-au^2$ or $au^2-t^2$, giving to $a$ any of the values up to 12. This method is, as we see, very easy and simple to use, but it seems subject to a kind of inconvenience: it sometimes gives more formulas than necessary to represent all the divisors of the numbers of the given form, so that it happens that some of these formulas become the same [thing], as we have seen it in the preceding examples. To remedy this, it would be necessary to have a general rule by which we could easily recognize the formulas which are mutually identical; it is this which we will examine, with all the generality to which the the material is susceptible. Since it is not demonstrated up until now that this identity of formulas 
cannot be among the divisors of the numbers of the form $t^2-au^2$, though the different cases of Section 18 provide no example, to leave nothing to desire on this subject, we consider equally the formulas of the one and the other type. \\

\begin{center}
\textsc{Problem III.}
\end{center}

22. Given the formula 
\[
py^2+2qyz+rz^2,
\]
in which $y$ and $z$ are variables and $p$, $q$, $r$ are positive or negative numbers, subject to the conditions that $pr-q^2 = a$ ($a$ being a given positive number) and $2a$ is neither $>p$ nor $>r$ (ignoring the signs of $p$, $q$, and $r$): to find if this formula may be transformed into another of the same type, which is subject to the same conditions. \\

Since the transformation must be analogous to the proposed [formula] it is evident that we cannot employ substitutions other than these:
\[
y\;=\; Ms+Nx, \;\;\; z\;=\;ms+nx,
\]
$s$ and $x$ being two new variables, and $M$, $N$, $m$, and $n$ arbitrary numbers. Indeed, these substitutions give a transformation of the form
\[
Ps^2+2Qsx+Rx^2,
\]
in which
\begin{eqnarray*}
P &=& pM^2+2qMm+rm^2, \\
Q &=& pMN+q(Mn+Nm)+rmn, \\
R &=& pN^2+2qNn+rn^2,
\end{eqnarray*}
and it will only remain to see if we may determine the numbers $M$, $N$, $m$, and $n$ for which $PR-Q^2 = a$, and [for which] $2Q$ is neither $>P$ nor $>R$. \\

To satisfy the first condition, I substitute the values of $P$, $Q$, and $R$ into the quantity $PR-Q^2$, and I cancel terms\footnote{Literally, ``by erasing what is destroyed''.} [to obtain] $PR-Q^2 = (pr-q^2)(Mn-Nm)^2$. But $pr-q^2=a$ (by hypothesis), so for $PR-Q^2$ to also be $=a$ we must have $(Mn-Nm)^2 = 1$; consequently, $Mn-Nm = \pm 1$. \\

With regard to the second condition, it is clear that it cannot happen unless $Q$ is simultaneously $<P$ and $<R$; so we suppose that $Q$ is indeed $<P$ and $<R$, and we will see what must follow. \\

Let $M$ be greater than $N$ (the reasoning being the same if $N$ were greater than $M$, merely taking $N$ in place of $M$) it is clear that we may make
\[
M \;=\; \mu N + M',
\]
and that we may take $\mu$ so that $M'$ is less than $N$, since we need only take for $\mu$ the quotient of the division of $M$ by $N$, and $M'$ will be the remainder. Furthermore, it is easy to see that we may always suppose that $\mu$ is not less than $2$; so if we found $\mu =1$ such that $M = N+M'$, we could make $M = 2N-(N-M')$, which is to say take $\mu = 2$ and $N-M'$ in place of $M'$. But if we suppose also---which is permitted---that
\[
m \;=\; \mu n + m',
\]
$m'$ being any number, and that we substitute these values of $M$ and $m$ into the expression for $Q$, it will become
\[
Q \;=\; \mu(pN^2+2qNn+rn^2)+pM'N+q(M'n+Nm')+rm'n,
\]
or by making the abbreviation
\[
Q' \;=\; pM'N + q(M'n+Nm')+rm'n
\]
we will have
\[
Q \;=\; \mu R+Q'.
\]
Now it must be that $Q < R$; since $\mu \ge 2$, it is clear that this condition cannot happen unless the two quantities $\mu R$ and $Q'$ have different signs 
and at the same time $Q'$ is greater than $R$ (ignoring the signs).
\\

Now we will have $y=(\mu N+M')s+Nx$ and $z=(\mu n+m')s+nx$; so that if we make
\[
x' \;=\; \mu s+x
\]
we will have
\[
y \;=\; M's+Nx', \;\;\;\; z\;=\; m's+nx',
\]
and the substitution of these values into the formula $py^2+2qyz+rz^2$ will give the new transformed [expression] 
\[
P's^2 + 2Q'sx'+R{x'}^2,
\]
with
\begin{eqnarray*}
P' &=& p{M'}^2+2qM'm', \\
Q' &=& pM'N+q(M'nNm')+rm'n, \\
R &=& pN^2+2qNn+rn^2,
\end{eqnarray*}
as done above. \\

Now because $Mn-Nm = \pm 1$ we will have $(\mu N+M')n-N(\mu n+m') = \pm 1$, and consequently 
\[
M'n - Nm' \;=\; \pm 1.
\]
We will also find that $P'R-{Q'}^2 = (pr-q^2)\times (M'n-Nm')^2$, and consequently
\[
P'R-{Q'}^2 \;=\; a.
\]
So then, since $a$ is positive and $Q' > R$, it must be that $P' < Q'$; so the preceding transformed [expression] will be such that $R < Q' < P'$. \\ 

In the same fashion, because $N>M'$ we may suppose
\[
N \;=\; \mu'M'+N',
\]
and take $\mu'$ to be not $<2$, and $N' < M'$, and then make
\begin{eqnarray*}
n &=& \mu'm'+n', \\
s' &=& \mu'x'+s,
\end{eqnarray*}
so that we have 
\[
y \;=\; M's'+N'x', \;\;\;\; z \;=\; m's'=n'x'.
\]
We will obtain, using similar operations and reasoning as before, this new transformed [expression]
\[
P'{s'}^2+2Q''s'x'+R'{x'}^2,
\]
from which we will have
\begin{eqnarray*}
P' &=& p{M'}^2 + 2qM'm' +r{m'}^2, \\
Q'' &=& pM'N' + q(M'n'+N'm') + rm'n', \\
R' &=& p{N'}^2 + 2qN'n' + r{n'}^2,
\end{eqnarray*}
and where we will also have
\begin{eqnarray*}
M'n' - N'm' &=& \pm 1, \\
Q' &=& \mu'P'+Q', \\
P'R' -{Q''}^2 &=& a,
\end{eqnarray*}
so that $Q''$ will be $>P'$ and $< R'$, ignoring the signs of $P$, $R'$, and $Q''$. \\

We can also find a third transformation, such as
\[
P''{s'}^2+2{Q'''}s'x''+R'{x''}^2,
\]
which will be subject to the same conditions of the preceding transformations, and so on. \\

I now consider that, as the numbers $M$, $N$, $M'$, $N'$, etc. form (ignoring their signs) a decreasing sequence, we will necessarily arrive at a term which will be $=0$. Suppose that $N'$ is this term, so that we have $N'=0$; then because $M'n'-N'm' = \pm 1$ we will have $M'n' = \pm 1$; thus $M'=\pm 1$ and $n'=\pm 1$, therefore
\begin{eqnarray*}
P' &=& p\pm 2qm'+r{m'}^2, \\
Q'' &=& \pm q\pm rm', \\
R' &=& r,
\end{eqnarray*}
the $\pm$ signs being arbitrary. \\

Now it must be (1) that we have $Q'' < R'$, ignoring the signs of these numbers; but $R'=r$ and $q<r$, because $2q$ is not $>r$ (by hypothesis); thus $Q''$ cannot be $< R' < r$ unless $m'$ is $=0$ our $=\pm 1$. (2) It must be that $Q''$ is at the same time $>P'$; but if $m'=0$ we have $Q''=\pm q$; and $P'=p$; so that because $2q$  is not $>p$ (by hypothesis), $Q''$ will always be $<P'$ instead of being greater; if $m'=\pm 1$ we will have $P'=p\pm 2q+r$, and $Q''=q\pm r$; but we suppose that $Q''$ will be $<r$; thus, for $Q''$ to be $>P'$ it would be necessary for $r$ to be $>p\pm 2q+r$, which cannot happen because $2q$ is never $>p$, and besides, $p$ and $r$ must have the same signs by virtue of the equation $pr-q^2$ equaling a positive number. \\

From this I conclude that it is impossible for the proposed formula to be transformed into another where the stated conditions take place; so that if we have several formulas where the same conditions are observed, we may be assured that the formulas are essentially different from each other, and that they cannot be reduced to an even smaller number. \\

\begin{center}
\textsc{Problem IV.}
\end{center}

23. \emph{Being given the formula
\[
py^2+2qyz-rz^2,
\]
in which $y$ and $z$ are variables and $p$, $q$, $r$ are positive (or negative) numbers determined by these conditions, that $pr+q^2 = a$ ($a$ being a given positive number), and that $2q$ is neither $>p$ nor $>r$, ignoring the signs of $p$, $q$, and $r$; to find if this formula can be transformed into another similar [one], where the same conditions are observed.} \\

Proceeding as in the preceding Problem, and for the same reason,
\[
y\;=\;Ms+Nx,\;\;\;\; z\;=\;ms+nx,
\]
we will have the transformed [expression] 
\[
Ps^2+2Qsx-Rx^2,
\]
in which
\begin{eqnarray*}
P &=& pM^2+2qMm-rm^2, \\
Q &=& pMN+q(Mn+Nm)-rmn, \\
R &=& rn^2-2qNn-pN^2,
\end{eqnarray*}
so the difficulty consists in determining (if it is possible) the numbers $M$, $N$, $m$, and $n$ for which we have $PR+Q^2 = a$, and at the same time neither $P$ nor $R$ being $<2Q$, ignoring the signs of $P$, $Q$, and $R$. \\

I note first that, by putting in place of $P$, $Q$, and $R$ their [respective] values, the quantity $PR+Q^2$ becomes $(pr+q^2)(Mn-Nm)^2=a(Mn-Nm)^2$; so necessarily we will have, as in the preceding Problem, $(Mn-Nm)^2 = 1$, and consequently
\[
Mn-Nm\;=\;\pm 1.
\]
Since $M$, $N$, $m$, and $n$ are assumed to be whole numbers, it is clear that this equation cannot exist unless the products $Mn$ and $Nm$ do not have the same signs; so then if $M$ and $N$ have the same signs, it must be that $m$ and $n$ do as well. \\

Now, since we may give variables $s$ and $x$ such signs as we want, it is evident that we may, without harm to the generality of the Problem, always take the numbers $M$ and $N$ to be positive; and then it must be that the numbers $m$ and $n$ have the same sign, which is to say both are positive, or both are negative. Then it will only be necessary to put $\pm m$ and $\pm n$ in place of $m$ and $n$, or, what amounts to the same thing, there will be no need to take the ambiguous sign $\pm$ in the quantity $q$, which is to say to take the value of this quantity in \emph{plus} and in \emph{minus}; through which we may regard the four numbers $M$, $N$, $m$, and $n$ as positive. \\

Now it is clear that if $2Q$ is neither $>P$ nor $>R$, as we have supposed, $Q^2$ will always be less than $PR$, so that $PR+Q^2$ cannot be equal to a positive number unless $PR$ is a positive number; from which it follows that it is necessary for $P$ and $R$ to have the same sign; and this condition suffices, as we will see, to find the numbers $M$, $N$, $m$, $n$. \\

[To show] this, I observe that because $pr+q^2 = a$, the quantity $P$ can be put in this form
\[
P \;=\; p\left(M+\frac{q+\sqrt{a}}{p}\,m\right)\left(M+\frac{q-\sqrt{a}}{p\,}m\right),
\]
and the quantity $R$ in this one
\[
R \;=\; -p\left(N+\frac{q+\sqrt{a}}{p}\,n\right)\left(N+\frac{q-\sqrt{a}}{p}\,n\right).
\]
Now, since $\sqrt{a}$ is greater than $q$, it is clear that the quantity $q+\sqrt{a}$ will always be positive, and the quantity $q-\sqrt{a}$ always negative; so that the two quantities $\frac{q+\sqrt{a}}{p}$ and $\frac{q-\sqrt{a}}{p}$ will necessarily have different signs. Letting then $\alpha$ be that of the two quantities which is positive and $-\beta$ the one which is negative ($\alpha$ and $\beta$ denoting positive numbers), we will have
\begin{eqnarray*}
P &=& p(M+\alpha m)(M-\beta m), \\
R &=& -p(N+\alpha n)(N-\beta n).
\end{eqnarray*}
From this we see that for the numbers $P$ and $R$ to have the same sign it must be that the factors $M-\beta m$ and $N-\beta n$ must have different signs, because the factors $M+\alpha m$ and $N+\alpha n$ are both positive. \\

That said, let $M>N$; we may make $M = \mu N+M'$ and take for $\mu$ a positive whole number for which $M'$ will be positive and less than $N$; because for that we need only divide $M$ by $N$ and make the quotient equal to $\mu$ and the remainder equal to $M'$. Let us do the same with $m=\mu n+m'$, $m'$ being an arbitrary number; and substituting these values into the equation $Mn-Nm=\pm 1$, we will have this
\[
M'n - Nm' \;=\; \pm 1,
\]
where we see that because $M'$, $N$, and $n$ are positive, it must be that $m'$ is also a positive number. \\

Now, the values of $y$ and of $z$ will become by the same substitutions 
\[
y \;=\; (\mu s+x)N+M's, \;\;\;\; z \;=\; (\mu s+x)n+m's,
\]
or else, by doing as above 
\begin{eqnarray*}
x' &=& \mu s+x, \\
y &=& M's+Nx', \\
z &=& m's+nx',
\end{eqnarray*}
and these values being substituted into the formula $py^2+2qyz-rz^2$, we will have the transformed [expression]
\[
P's^2+2Q'sx'-R{x'}^2
\]
where
\begin{eqnarray*}
P' &=& p{M'}^2+2qM'm'-r{m'}^2, \\
Q' &=& pM'N+q(M'n+Nm')-rm'n, \\
 R &=& rn^2-qNn-pN^2.
\end{eqnarray*}
And I say that the numbers $P'$ and $R$ will necessarily have the same signs; now we will have
\begin{eqnarray*}
P' &=& p(M'+\alpha m')(M'-\beta m'), \\
R &=& -p(N+\alpha n)(N-\beta n);
\end{eqnarray*}
so $M-\beta m = \mu(N-\beta n) + M'-\beta m'$. Thus, since $\mu$ is a positive number, and [since] $M-\beta m$ and $N-\beta n$ have different signs, for this equation to exist it must be that the quantities $M-\beta m$ and $M'-\beta m'$ have the same signs; and consequently that $N-\beta n$ and $M'-\beta m'$ have different signs; but $N+\alpha n$ and $M'+\alpha m'$ are positive quantities, $N$, $n$, $m'$, and $\alpha$ being positive numbers; thus the two numbers $P'$ and $R$ necessarily have the same sign. \\

Similarly, since $N>M'$, we may suppose $N=\mu' M'+N'$ and take $\mu'$ positive so that $N'$ is also positive and less than $M'$; and making $n=\mu' m'+n'$ we will have (substituting these values into the equation $M'n-Nm' = \pm 1$)
\[
M'n' - N'm' \;=\; \pm 1,
\]
so that $m'$ will also be necessarily positive.
\\

Next, if we make 
\[
s' \;=\; \mu'x'+s,
\]
we will have 
\[
y\;=\;M's'+N'x',\;\;\;\; z\;=\;m's'+n'x',
\]
and substituting these values into the formula $py^2+2qyz-rz^2$, we will have this other transformed [expression]
\[
P'{s'}^2+2Q''s'x'-R'{x'}^2,
\]
where 
\begin{eqnarray*}
P' &=& p{M'}^2 + 2qM'm' - r{m'}^2, \\
Q'' &=& pM'N' + q(M'n'+N'm')-rm'n', \\
R' &=& r{n'}^2-qN'n'-p{N'}^2.
\end{eqnarray*}
And we will prove, as we have done so far, that the numbers  $P'$ and $R'$ will have the same signs. \\

Similarly, we can find a third transformation such that
\[
P''{s'}^2 + 2Q'''s'x'' - R'{x''}^2,
\]
in which 
\[
x'' \;=\; \mu''s'+x',
\]
and where $P''$ and $R'$ will have the same signs, and so on. \\

Now, since the numbers $M$, $N$, $M'$, $N'$ etc. form a decreasing sequence of whole numbers, it is clear that we must necessarily reach a term which is zero. Suppose then, for example, that we have $N'=0$, and because $M'n'-N'm' = \pm 1$, we will have $M'n'=1$ (so, because the numbers $M'$ and $n'$ are both positive, it is evident that we must take the positive sign in this case), thus $M'=1$ and $n'=1$; so that we will have in this case $y=s'$, $z=m's'+x'$. \\ 

Therefore I conclude that, to transform the proposed formula
\[
py^2+2qyz-rz^2
\]
into this one,
\[
Ps^2+2Qsx-Rx^2,
\]
in which we have $PR+Q^2=pr+q^2=a$, and where $P$ and $R$ have the same signs, we must do the following substitutions
\begin{eqnarray*}
z &=& m'y + x', \\
y &=& \mu'x' + s, \\
x &=& \mu s + x,
\end{eqnarray*}
and take the numbers $m'$, $\mu'$, and $\mu$ to be positive, and such that in the resulting transformations
\begin{eqnarray*}
& & P'y^2 + 2Q''yx' - R'{x'}^2, \\
& & P's^2 + 2Q'sx' - R{x'}^2, \\
& & Ps^2 + 2Qsx - Rx^2,
\end{eqnarray*}
the coefficients $R'$, $P'$, $R$, and $P$ have the same signs. \\

Let us see how we can fulfill these conditions. \\

First doing the substitution of $m'y+x$ in place of $z$, we will have the first transformed [expression], where
\begin{eqnarray*}
R' &=& r, \\
Q'' &=& q-rm', \\
P' &=& p+2qm' - r{m'}^2 \;=\; \frac{a-{Q''}^2}{R'}.
\end{eqnarray*}
Now, $P' = -r \left(m'+\frac{\sqrt{a}-q}{r}\right)\left(m'-\frac{\sqrt{a}+q}{r}\right)$. Thus, in order for $P'$ and $R'$ to have the same signs, the factors $m'+\frac{\sqrt{a}-q}{r}$ and $m'-\frac{\sqrt{a}+q}{r}$ must have different signs, but because $\sqrt{a} > q$, it is clear that $\sqrt{a}\pm q$ will always be a positive number; thus, if $r$ is positive, $m'+\frac{\sqrt{a}-q}{r}$ will always be positive, and it must be that $m'-\frac{\sqrt{a}+q}{r}$ is negative, and consequently that $m' < \frac{\sqrt{a}+q}{r}$. If, to the contrary, $r$ is negative, $m'-\frac{\sqrt{a}+q}{r}$ will be positive and it must be that $m'+\frac{\sqrt{a}-q}{r}$ is negative; thus, $m' < \frac{\sqrt{a}-q}{-r}$. \\

Next substitute $\mu'x'+s$ in place of $y$ and we will have the second transformed [expression] in which
\begin{eqnarray*}
Q' &=& Q'' + P'\mu' \\
R &=& R' - 2Q''\mu' -P'{\mu'}^2 \;=\; \frac{a-{Q'}^2}{P'}.
\end{eqnarray*}
I observe that $R = -P' \left(\mu'+\frac{\sqrt{a}+Q''}{P'}\right)\left(\mu'-\frac{\sqrt{a}-Q''}{P'}\right)$, so that for $R$ and $P'$ to have the same signs it must be that the two factors $\mu'+\frac{\sqrt{a}+Q''}{P'}$ and $\mu'-\frac{\sqrt{a}-Q''}{P'}$ have different signs; so since $P'R' = a-{Q''}^2$ ($P'$ and $R'$ having the same signs) it follows that ${Q''}^2 < a$, and consequently $Q'' < \sqrt{a}$, so that $\sqrt{a}\pm Q''$ will always be a positive number; thus, if $P'$ is positive, $\mu'+\frac{\sqrt{a}+Q''}{P'}$ will be positive, and it must be that $\mu'-\frac{\sqrt{a}-Q''}{P'}$ is negative; thus $\mu' < \frac{\sqrt{a}-Q''}{P'}$; but for $\mu'$ to be a whole number it must be that $\frac{\sqrt{a}-Q''}{P'} > 1$; thus $P'<\sqrt{a}-Q''$; therefore, because $P'R'=a-{Q''}^2 = (\sqrt{a}+Q'')(\sqrt{a}-Q'')$, it must be that $R' > \sqrt{a}+Q''$, which is to say $r>\sqrt{a}+q-rm'$; and consequently $(m'+1)r > \sqrt{a}+q$ and then $m' > \frac{\sqrt{a}+q}{r}-1$. So $P'$ must be positive when $r$ is positive, in which case we have already found $m' < \frac{\sqrt{a}+q}{r}$; therefore we will have in this case
\begin{eqnarray*}
m' &<& \frac{\sqrt{a}+q}{r} \text{ and } > \frac{\sqrt{a}+q}{r}-1, \\
\mu' &<& \frac{\sqrt{a}-Q''}{P'}.
\end{eqnarray*}
We will find the same for the case when $r$ is negative
\begin{eqnarray*}
m' &<& \frac{\sqrt{a}-q}{-r} \text{ and } > \frac{\sqrt{a}-q}{-r}-1, \\
\mu' &<& \frac{\sqrt{a}+Q''}{P'}.
\end{eqnarray*}
From this we see that the number $m'$, needing to be whole, will be necessarily determined, since the two limits between which it must be found only differ by a unit. \\

Finally we will substitute $\mu s+x$ in place of $x'$, and we will have the third transformed [expression] in which
\begin{eqnarray*}
Q &=& Q' - R\mu, \\
R &=& P' + 2Q'\mu - R\mu^2 \;=\; \frac{a-Q^2}{R}.
\end{eqnarray*}
And paying attention to $p = $ (because $RP' = a-{Q'}^2$) we will prove, as above, that in the case where $r$ is positive, we will have 
\begin{eqnarray*}
\mu' &<& \frac{\sqrt{a}-Q''}{P'} \text{ and } > \frac{\sqrt{a}-Q''}{P'}-1, \\
\mu &<& \frac{\sqrt{a}+Q'}{R},
\end{eqnarray*}
and in the case where $r$ is negative,
\begin{eqnarray*}
\mu' &<& \frac{\sqrt{a}+Q''}{P'} \text{ and } > \frac{\sqrt{a}+Q''}{P'}-1, \\
\mu &<& \frac{\sqrt{a}-Q'}{R}.
\end{eqnarray*}
So then the number $\mu'$ will also be determined, and the only indeterminate number will be $\mu$. \\

Now if we want further that $2Q$ is neither $>P$ nor $>R$, as in the conditions of the required Problem, we must first determine $\mu$ so that $Q = Q'-\mu R$ is not $>\frac{R}{2}$, ignoring the signs of $Q$ and $R$; and it is clear that, taking $\mu$ to be a positive whole number, there is not a single value of $\mu$ which can satisfy this condition; so then the number $\mu$ will be completely determined by this method. So it only remains to see if $Q$ is also $< \frac{P}{2}$, in which case the transformed [expression] $Ps^2+2Qsx-Rx^2$ will have the required conditions. \\

From this comment, we see how the proposed question may be answered without any guesswork; and here is the method which must be followed to this end. \\

24. \emph{Method to transform the formula 
\[
py^2+2qyz-rz^2
\]
in which we have $pr+q^2=a$ ($a$ being a given, positive whole number) and where $2q$ is neither $>p$ nor $>r$ (ignoring the signs of $p$, $q$, $r$), into other similar formulas and subject to the same conditions.} \\

To better preserve the analogy 
in our formulas, we change first the letters $z$ and $p$ to $y'$ and $r'$, so that our formula will become 
\[
r'y^2 + 2qyy'-r{y'}^2,
\]
where $rr'+q^2 = a$, and $q$ is neither $>\frac{r}{2}$ nor $>\frac{r'}{2}$. \\

Now, since $r$ and $r'$ must have the same signs by virtue of the equation $r'r+q^2=a$, we will suppose first that they are both positive; but $q$ will have to be taken successively as positive and negative. \\

That said, we will have \\

$1^\circ.\; y=m'y'+y''$, which will give this first transformed [expression]
\[
r'{y''}^2+2q'y''y'-r''{y'}^2,
\]
where we will have 
\begin{eqnarray*}
q' &=& q+r'm', \\
r'' &=& r-2qm'-r'{m'}^2 \;=\; \frac{a-{q'}^2}{r'}.
\end{eqnarray*}
We will take $m'$, if possible, to be a positive whole number, such that $q+r'm'$ is not $>\frac{r'}{2}$; then we will see if $r'' > q'$ or not; and in this last case the transformed [expression thus] found will have the required conditions. \\

$2^\circ$. We will determine $m'$ so that 
\[
m' < \frac{\sqrt{a}-q}{r'} \text{ and } > \frac{\sqrt{a}-q}{r'}-1.
\]
Next we will take
\[
y' \;=\; m''y''+y''',
\]
which will give this second transformed [expression]
\[
r'''{y''}^2 + 2q''y''y'''-r''{y'''}^2,
\]
by making
\begin{eqnarray*}
q'' &=& q'-r''m'', \\
r''' &=& r'+2q'm''-r''{m''}^2 = \frac{a-{q''}^2}{r''}.
\end{eqnarray*}
We will take $m''$ to be a positive whole [number] such that $q'-r''m''$ is not $>\frac{r''}{2}$; and if at the same time $q''$ does not exceed $\frac{r'''}{2}$ the preceding transformed [equation] will have the required conditions. \\

$3^\circ$. We will determine $m''$ so that 
\[
m'' < \frac{\sqrt{a}+q'}{r''} \text{ [and] } >\frac{\sqrt{a}+q'}{r''} - [1].
\]
Next we will take
\[
y'' \;=\; m'''y''' + y^{IV},
\]
and we will have this third transformed [expression]
\[
r'''{y^{IV}}^2 + 2q'''y^{IV}y''' - r^{IV}{y'''}^2,
\]
in which
\begin{eqnarray*}
q''' &=& q'' + r'''m''', \\
r^{IV} &=& r''-2q''m'''-r'''{m'''}^2 \;=\; \frac{a-q''}{r'''}.
\end{eqnarray*}
We will take for $m'''$ a positive whole number for which $q''+r'''m'''$ is not $>\frac{r'''}{2}$, and if the value of $q'''$ is not simultaneously $>\frac{r^{IV}}{2}$, we will be assured that the transformed [expression thus] found will have the required conditions. \\

$4^\circ$. We will determine $m'''$ so that 
\[
m''' < \frac{\sqrt{a}-q''}{r'''} \text{ and } >\frac{\sqrt{a}-q''}{r'''} - 1.
\]
Next we will take
\[
y''' \;=\; m^{IV}y^{IV} + y^{V},
\]
which will give the fourth transformed [expression]
\[
r^{V}{y^{IV}}^2 + 2q^{IV}y^{IV}y^{V} - r^{IV}{y^V}^2,
\]
where
\begin{eqnarray*}
q^{IV} &=& q''' + r^{IV}m^{IV}, \\
r^{V} &=& r'''-2q'''m^{IV}-r^{IV}{m^{IV}}^2 \;=\; \frac{a-{q^{IV}}^2}{r^{IV}}.
\end{eqnarray*}
We will take $m^{IV}$ so that $q''-r^{IV}m^{IV}$ is not $>\frac{r^{IV}}{2}$, and if $q^{IV}$ is not simultaneously $>\frac{r^V}{2}$ the transformed [expression] will have the required conditions. \\

$5^\circ$. We will determine $m^{IV}$, etc. \\

In this fashion we will find successively all the transformations of the proposed formula in which the prescribed conditions can hold; and it is clear that the number of different transformations will be necessarily limited; now we have seen in Prob. II that we can only have a limited number of different formulas where the same conditions are observed. \\

But to get all the possible different transformations of the same formula, it will be necessary to make a double calculation taking the value of $q$ successively as positive and negative. \\

If the numbers $r$ and $r'$, instead of being two positive numbers, as we have supposed, are two negatives, it would only be necessary to change the signs of these numbers as well as the number $q$, which is to say that we would take the formula $r'y^2+2qyy'-r'y^2$ to be negative; and then we would at the same time change all the signs of the transformed [expressions] that we would have found. Or, which is even simpler, we will write $-r$ in place of $r'$, $-r'$ in place of $r$, and $y'$ in place of $y$, which will give the formula $-r{y'}^2+2qyy'+r'y^2$ where $r$ and $r'$ will be positive numbers. 

\begin{center}
\textsc{Corollary.}
\end{center}

25. It follows from the analysis of the preceding Problem that the numbers $r$, $r'$, $r''$, etc. will all have the same signs and satisfy $r''+q^2 = a = r'r''+{q'}^2 = r''r'''+{q''}^2 = $ etc., so each of these numbers will be less than the given number $a$. Consequently, in continuing the series $r$, $r'$, $r''$ etc. it will be necessary that the same number repeats several times and even that the same pair of successive numbers also repeats; thus in continuing the calculation, following the preceding \emph{method}, we will necessarily recover a transformed [expression] identical to one of those which we already had. It is this which we will recognize easily when we will find, for example, $q^{\mu+\nu} = q^{\mu}$ and $r^{\mu+\nu+1} = r^{\mu + 1}$, and that $\nu$ will be an even number; so it will be useless to pursue the calculation any further, because the following transformations will be the same as those which we have already found. \\ 

Therefore, as soon as we have found from Problem II all the different formulas $py^2\pm 2qyz - rz^2$ which may represent the divisors of the numbers of the form $t^2-au^2$, we may reduce them to the smallest possible number, excluding those which are not the transformed [expressions] of some one of these formulas. So, since the formula $y^2-az^2$ is always one of the divisors of $t^2-au^2$ (by making $q=0$ and $p=1$, $r=a$) we will begin by searching all the transformed [expressions] of this same formula, where the prescribed properties will occur, and as these transformed [expressions] are found necessarily among the other formulas of the divisors of $t^2-au^2$ we may first reject those which are identical among themselves. Next we will do the same operation on the formulas which remain; and after those have been covered, rejecting those which are found to be identical among themselves, we will be sure that those remaining will be necessary to represent all the possible divisors of the numbers of the given form. \\ 

Moreover, it will happen most often that the transformations of the formula $y^2-az^2$ will contain all the other formulas of the divisors of $t^2-au^2$, especially when $a$ is a prime number; but it would be wrong to make a general rule, so we will give examples where they find themselves at fault, which will serve simultaneously to show the utility and importance of the methods which we just gave. \\

\begin{center}
\textsc{Examples.}
\end{center}

26. Let us propose the formula $y^2-2z^2$; so $r'=1$, $q=0$, $r=2=a$; so we will have $q'=m'$, $r''=\frac{2-{q'}^2}{1}$; so it is clear that we cannot make $q'<r'<1$; so we will pass to a second transformed [expression]. \\

For that we will thus take $m' < \frac{\sqrt{2}}{1}y$ [and] $> \frac{\sqrt{2}}{1} - 1$, which is to say $m'=1$, which will give $q'=1$, $r''=\frac{2-1}{1} - 1$; next we will have $q''=q'-r''m''=1-m''$ and $\frac{2-{q''}^2}{1}$; so, for $q''$ to not be $>\frac{r''}{2} > \frac{1}{2}$, it must be that $m''=1$, which gives $q''=0$ and $r'''=2$; so that as $q''$ is at the same time not $\frac{r''}{2}$ we will have the transformed [expression] $r'''{y''}^2+2q''y''y''' - r''{y'''}^2$, which is to say 
$2{y''}^2-{y'''}^2$, which will have the required conditions. Now this transformed [expression] is similar to the formula $2z^2-y^2$, so that the two formulas $y^2-2z^2$ and $2z^2-y^2$, which our general method gives for the divisors of the numbers of the form $t^2-2u^2$, return to the same 
as we have already remarked (Article 20, No. II). \\

We will find, likewise, that the two formulas $y^2-5z^2$ and $5z^2-y^2$ return to the same, 
as we observed it in the Article cited (No. V). \\

To give another example, consider the case of the No. VII of the same Article, where we have found that the formulas of the divisors of $t^2-7u^2$ were $y^2-7z^2$, $2y^2\pm 2yz-3z^2$, $7z^2-y^2$, and $3z^2\pm 2yz -2y^2$. \\

So first let $r'=1$, $q=0$, and $r=7=a$; we will have $q'=m'$ and $r''=\frac{7-{q'}^2}{1}$, from which we see that $q'$ cannot be $< \frac{r'}{2} < \frac{1}{2}$. \\

Second we will take $m' < \frac{\sqrt{7}}{1}$ and $>\frac{\sqrt{7}}{1}-1$, so $m'=2$ and $q'=2$ and $r''=3$; from here we will have $q''=2-3m''$, $r'''=\frac{7-{q''}^2}{3}$, and since $q''$ is not $>\frac{r''}{2} > \frac{3}{2}$ it must be that $m''=1$, which will give $q''=-1$, and $r'''=2$. So since $q''$ is at the same time not $>\frac{r'''}{2}$, the transformed [expression is] $r'''{y''}^2+2q''y''y'''-r''{y'''}^2$, which is to say [that] $2{y''}^2-2y''y'''-3{y'''}^2$ will have the required conditions. \\

Third we will take $m'' < \frac{\sqrt{7}+2}{3}$ and $>\frac{\sqrt{7}+2}{3}-1$, so $m''=1$; from which $q'' = -1$ and $r'''=2$. Then we will have $q''' = -1+2m'''$ and $r^{IV} = \frac{7-{q'''}^2}{2}$; thus for $q'''$ to not be $>\frac{r'''}{2}$ we must take $m'''=1$. This will give $q'''=1$ and $r^{IV}=3$; from that we will have the new transformed [expression] $2{y^{IV}}^2+2y^{IV}y'''-3{y'''}^2$, which will also have the required conditions. \\

Fourth we will take $m''' < \frac{\sqrt{7}+1}{2}$ and $>\frac{\sqrt{7}+1}{2}-1$; thus $m'''=1$, from which $q'''=1$, $r^{IV}=3$; next we will have $q^{IV} = 1-3m^{IV}$ and $r^V = \frac{7-{q^{IV}}^2}{3}$, where we see that we cannot take $m^{IV}$ so that $q^{IV}$ does not become $>\frac{r^{IV}}{2}$. \\

Fifth we will take $m^{IV} < \frac{\sqrt{7}+1}{3}$ and $>\frac{\sqrt{7}+1}{3}-1$; thus $m^{IV} = 1$, and from here $q^{IV} = -2$ and $r^V = 1$. Next we will have $q^V = -2+m^{V}$, $r^{VI} = \frac{7-{q^V}^2}{1}$; thus for $q^V$ to not be $>\frac{r^{V}}{2}$ we will be $m^V = 2$, which will give $q^V = 0$ and $r^{VI} = 8$; so that we will have the transformed [expression] ${y^{VI}}^2-7{y^V}^2$, which will have the prescribed conditions. \\

Sixth we will take $m^{V} < \frac{\sqrt{7}+2}{1}$ and $>\frac{\sqrt{7}+2}{1}-1$; thus $m^V = 4$, consequently $q^V = 2$ and $r^{VI} = 3$. Here I observe, without going further, that these values of $q^V$ and $r^{VI}$ are the same as $q'$ and $r''$ (No. 2); thus, since the difference of the exponents of $q$ is even, it follows that the transformed [expressions] which we could find by continuing the calculation would be the same as those which we have already found above (Article 25). \\

So the formula $y^2-7{y'}^2$ cannot furnish other transformations which have these two prescribed conditions: $2{y'''}^2-2y'''y'''-3{y'''}^2$ 
and $2{y^{IV}}^2 +2{y^{IV}}y'''-3{y'''}^2$; from which we see that the formulas $y^2-7z^2$ and $2y^2\pm 2yz-3z^2$ revert to the same thing, since the formulas $7z^2-y^2$ and $3z^2\mp 2yz-2y^2$ are only the negatives of those. But the two formulas $y^2-7z^2$ and $7z^2-y^2$ cannot reduce from the one to the other, as took place in the formulas $y^2-5z^2$ and $5z^2-y^2$ from the preceding example. \\

27. To further develop the application of our methods from Problems II and IV, we will seek here the formulas of the divisors of the numbers of the form $t^2-79u^2$ or $79u^2-t^2$. \\ 

Here we will then have $a=79$; so it must be that $q$ is not $>\sqrt{\frac{79}{3}} > 3$, so that we may only do $q=0$, $1$, $2$, $3$. Making $q=0$ we will have $pr=79$, so $p=1$ and $r=79$; making $q=1$ we will have $pr=78$; so either $p=2$ and $r=39$, or $p=3$ and $r=26$, or $p=6$ and $r=13$; making $q=2$ we will have $pr=75$, so $p=5$ and $r=15$; lastly making $q=3$, we will have $pr=70$, so $p=7$ and $r=10$. \\

So for the divisors in question, we will have the following formulas, $y^2-79z^2$, $2y^2\pm 2yz-39z^2$, $3y^2\pm 2yz-26z^2$, $6y^2\pm 2yz-13z^2$, $5y^2\pm 4yz -15z^2$, $7y^2\pm 6yz -10z^2$, and their inverses $79y^2-z^2$, $39z^2\mp 2yz-2y^2$, $26z^2\mp 2yz-3y^2$, $13z^2\mp 2yz-6y^2$, $15z^2\mp 4yz-5y^2$, $10z^2\mp 6yz-7y^2$, which makes 12 formulas in total; but we must now sort them, and reject those which are identical to one another. \\ 

Consider first the formula $y^2-79z^2$, or $y^2-79{y'}^2$, and we will have (1) $r'=1$, $q=0$, and $r=79=a$, so $q'=m'$ and $r''=\frac{79-{q'}^2}{1}$; now $q'$ is always $>\frac{r'}{2}$, unless we make $m'=0$, which would give no new formula. \\

So we will have (2) $m' < \frac{\sqrt{79}}{1}$ and $>\frac{\sqrt{79}}{1}-1$, thus $m'=8$, consequently $q'=8$ and $r''=15$; next we will have $q'' = 8-15m''$, $r'''=\frac{79-{q''}^2}{15}$, from where it will be $m''=1$ so that $q''$ is not $>\frac{r''}{2}$. We will thus have $q''=-7$ and $r'''=2$; but since $q''$ would be $>\frac{r'''}{2}$ these values give no suitable transformation. \\

Thus we will have (3) $m'' < \frac{\sqrt{79}+8}{15}$ and $>\frac{\sqrt{79}+8}{15}-1$; so $m''=1$ and from that $q''=-7$, $r''=2$; then we will have $q'''=-7+2m'''$ and $r^{IV} = \frac{79-{q'''}^2}{2}$. We take $m'''=3$ or $=4$ so that $q'''=\pm 1$ is not $>\frac{r'''}{2}$, and $r^{IV}$ will become $=39$; so that we will have this transformed [expression], which will have all the prescribed conditions: $2{y^{IV}}^2\pm 2y^{IV}y'''-39{y'''}^2$. \\

Continuing the calculation we will have (4) $m''' < \frac{\sqrt{79}+7}{2}$ and $>\frac{\sqrt{79}+7}{2}-1$, which is to say $m'''=7$; from which $q'''=7$ and $r^{IV} = 15$. Next, we will have $q^{IV}=7-15m^{IV}$ and $r^{V} = \frac{79-{q^{IV}}^2}{15}$; and we will take $m^{IV} = 1$ so that $q^{IV}$ is not $>\frac{r^{IV}}{2}$. So we will have $q^{IV} = -8$ and $r^{V}=1$; but since $q^{IV} > \frac{r^V}{2}$ we will reject these values as useless. \\

We will then have (5) $m^{IV} < \frac{\sqrt{79}+7}{15}$ [and] $\frac{\sqrt{79}+7}{15}-1$, so $m^{IV} = 1$, consequently $q^{IV} = -8$, $r^{V}=1$; after which we will suppose $q^V=-8+m^V$ and $r^{VI}=\frac{79-{q^V}^2}{1}$, and we will take $m^V = 8$ so that $q^V=0$ and $r^{VI} = 79$, which will give the transformed [expression] ${y^{VI}}^2-79{y^V}^2$ which is entirely similar to the first formula ${y'}^2-79y^2$. \\

I do (6) $m^V < \frac{\sqrt{79}+8}{1}$ and $>\frac{\sqrt{79}+8}{1}-1$, knowing $m^V =16$, which gives $q^V=8$ and $r^{VI}=15$. So I remark that these values of $q^V$ and $r^{VI}$ are the same as $q'$ and $r''$ from (2); so that since the difference of the exponents of $q$ is even, we will recover the same transformation which we have already had; from which it follows that the formula ${y'}^2-79y^2$ cannot be changed into any other [formula] than this one: $2y^{IV}\pm 2y^{IV}y'''-39{y'''}^2$. So, among all the formulas found for the divisors of $t^2-79u^2$ there are only these two: $y^2-79z^2$ and $2y^2\pm 2yz-39z^2$, which are identical, to which we must add their inverses $79z^2-y^2$ and $39z^2\mp 2yz-2y^2$, which are also identical. \\

Now consider the formula $3y^2\pm 2yz-26z^2$, knowing $3y^2\pm 2yy' - 26{y'}^2$, and we will have (1) $r'=3$, $q=1$, and $r=26$, $a$ always being $=79$. So we will suppose $q'=1+3m'$, $r''=\frac{79-{q'}^2}{3}$, and since we cannot take $m'$ so that $q'$ is not $>\frac{r''}{2}$, we will pass to another transformation. \\

We will then have (2) $m' < \frac{\sqrt{79}-1}{3}$ and $>\frac{\sqrt{79}-1}{3}-1$; so $m'=2$, $q'=7$, and $r''=10$; then we will suppose $q''=7-10m''$, $r'''=\frac{79-{q''}^2}{10}$; then we will take $m''=1$ to have $q''=-3 < \frac{10}{2}$, and we will have $r'''=7 > 2q''$. So we will have the transformed [expression] $7{y''}^2-6y''y'''-10{y'''}^2$, which will have the required conditions. \\

Let (3) $m'' < \frac{\sqrt{79}+7}{10}$ and $>\frac{\sqrt{79}+7}{10}-1$; thus $m''=1$, $q''=-3$, and $r'''=7$. Next let it be supposed that $q'''=-3+7m'''$ and $r^{IV} = \frac{79-{q'''}^2}{7}$, and since we cannot take $m'''$ so that $q'''$ is not $>\frac{r'''}{2}$, we will pass to the following transformation. \\

Thus we will have (4) $m''' < \frac{\sqrt{79}+3}{7}$ and $>\frac{\sqrt{79}+3}{7}-1$, which is to say $m'''=1$, and we will have $q'''=4$, $r^{IV}=9$. Next we will suppose $q^{IV} = 4-9m^{IV}$, $r^{V}=\frac{79-{q^{IV}}^2}{9}$, so we cannot take $m^{IV}$ so that $q^{IV}$ is not $>\frac{r^{IV}}{2}$; thus \footnote{The remainder of this sentence is missing in the original manuscript.} \\

We will have (5) $m^{IV} < \frac{\sqrt{79}+4}{9}$ and $>\frac{\sqrt{79}+4}{9}-1$, which is to say $m^{IV}=1$; thus $q^{IV} = -5$, and $r^V=6$; after which it will be $q^V=-5+6m^V$ and $r^{VI} = \frac{79-{q^V}^2}{6}$. Here we may take $m^V=1$, which gives $q^V=1$ and $r^{VI} = 13$, values which have the required conditions; so that we will have the transformed [expression] $6{y^{VI}}^2 + 2y^{VI}y^V - 13{y^V}^2$. \\

Now let (6) $m^V < \frac{\sqrt{79}+5}{6}$ and $>\frac{\sqrt{79}+5}{6}-1$, thus $m^V=2$, $q^V=7$, and $r^{VI}=5$; next let it be supposed $q^{VI} = 7-5m^{VI}$ and $r^{VII} = \frac{79-{q^{VI}}^2}{5}$, and it is clear that taking $m^{VI}=1$ we will have $q^{VI} < \frac{r^{VI}}{2}$. We will then have $q^{VI}=2$ and $r^{VII}=15$; so that the transformed [expression] $15{y^{IV}}^2+4y^{IV}y^V-5{y^V}^2$ will have the required conditions. \\

We take (7) $m^{VI} < \frac{\sqrt{79}+7}{5}$ and $>\frac{\sqrt{79}+7}{5}-1$, thus $m^{VI} = 3$ and $q^{VI}=-8$, $r^{VII}=3$; next let it be supposed $q^{VII} = -8+3m^{VII}$ and $r^{VIII} = \frac{79-{q^{VII}}^2}{3}$, and taking $m^{VII} =3$ we will have $q^{VII}=1 < \frac{r'''}{2}$ and $r^{VIII}=26>2q^{VII}$, which will give the transformed [expression] $3{y^{VIII}}^2+2y^{VIII}y^{VII} - 26{y^{VII}}^2$, which is similar to the proposed [expression]. \\

Next we will have (8) $m^{VII} < \frac{\sqrt{79}+8}{3}$ and $>\frac{\sqrt{79}+8}{3}-1$; that is to say $m^{VII} =5$, and consequently $q^{VII} = 7$ and $r^{VIII} = 10$; values which are the same as $q'$ and $r''$; so that the same transformed [expressions] which we have already found would return if we continued the calculation. \\

Now we represent the same values of $r'$ and $r$ from (1), but instead of supposing $q=1$, we do $q=-1$; so $q'=-1+3m'$ and $r''=\frac{79-{q'}^2}{3}$; so since we would not need to determine the $m'$ for which $q'$ becomes $<\frac{r'}{2}$, we must pass immediately to another transformation. \\

Thus we will have (2) $m' < \frac{\sqrt{79}+1}{3}$ and $>\frac{\sqrt{79}+1}{3} - 1$, so $m'=3$ and $q'=8$, $r''=5$; next we will suppose $q''=8-5m''$, $r'''=\frac{79-{q''}^2}{5}$; and it is clear that taking $m''=2$, $q''$ will not be $>\frac{r''}{2}$; so we will have $q''=-2$ and $r'''=15$, so that it will result in the transformation $15{y''}^2-4y''y'''-5{y'''}^2$, which has, as we see, the required conditions. \\ 

We will have (3) $m < \frac{\sqrt{79}+8}{5}$ and $>\frac{\sqrt{79}+8}{5} - 1$, which is to say $m''=3$; from where $q''=-7$, $r'''=6$. Next we will suppose $q'''=-7+6m'''$, $r^{IV} = \frac{79-{q'''}^2}{6}$; and we will take $m'''=1$ to have $q'''=-1$ and $r^{IV} = 13$, which will give the transformed [expression] $6{y^{IV}}^2-2y^{IV}y'''-13{y'''}^2$, which has the required conditions. \\

Let (4) $m''' < \frac{\sqrt{79}+7}{6}$ and $>\frac{\sqrt{79}+7}{6}-1$; so $m'''=2$ and $q'''=5$, $r^{IV} = 9$; next we will suppose $q^{IV} = 5-9m^{IV}$ and $r^{V} = \frac{79-{q^{IV}}^2}{9}$; next we will suppose $q^{IV} = 5-9m^{IV}$ and $r^V = \frac{79-{q^{IV}}^2}{9}$; and we may take $m^{IV} = 1$, which will give $q^{IV} = -4 < \frac{r^{IV}}{2}$. But then we will have $r^V = 7 < 2q^{IV}$, so that these values are not suitable. \\

So let (5)  $m^{IV} < \frac{\sqrt{79}+5}{9}$ and $>\frac{\sqrt{79}+5}{9}-1$, so $m^{IV} = 1$ and $q^{IV} = -4$, $r^V = 7$; then it will be $q^V = -4+7m^V$ and $r^{VI} = \frac{79-{q^V}^2}{7}$ and we may take $m^V = 1$, which will give $q^V=3<\frac{r^V}{2}$ and $r^{VI}=10>2q^V$. So we will have this transformation $7{y^{VI}}^2 +6y^{VI}y^V-10{y^V}^2$. \\

Let (6) $m^V < \frac{\sqrt{79}+4}{7}$ [and] $>\frac{\sqrt{79}+4}{7}-1$; so $m^V = 1$, and $q^V=3$, $r^{VI} = 10$; we then do $q^{VI}=3-10m^{VI}$ and $r^{VII} = \frac{79-{q^{VI}}^2}{10}$, and since we cannot take $m^{VI}$ so that $q^{VI}$ is not $>\frac{r^{VI}}{2}$ we will pass first to the following transformation. \\

So let (7) $m^{VI} < \frac{\sqrt{79}+3}{10}$ and $\frac{\sqrt{79}+3}{10}-1$, so $m^{VI}=1$ and $q^{VI}=-7$, $r^{VII}=3$; we then do $q^{VII} = -7+3m^{VII}$ and $r^{VIII} = \frac{79-{q^{VII}}^2}{3}$, and taking $m^{VII}=2$ we will have $q^{VII} = -1 < \frac{r^{VII}}{2}$, and $r^{VIII}=26>2q^{VII}$; thus we will have the transformed [expression] $3{y^{VIII}}^2-2{y^{VIII}}y^{VII}-26{y^{VII}}^2$, which is analogous to the proposed [formula]. \\

Finally let (8) $m^{VII} < \frac{\sqrt{79}+7}{3}$ and $>\frac{\sqrt{79}+7}{3}-1$, so $m^{VII}=5$ and $q^{VII}=8$, $r^{VIII}=5$, values which are the same as those of $q'$ and $r''$ in (2) above; so the operation will be terminated. \\

We see therefore that the formula $3y^2\pm 2yy'-26{y'}^2$ could only provide these transformations: $7{y''}^2-6y''y'''-10{y'''}^2$, $6{y^{V}}^2+2y^{VI}y^V-13{y^V}^2$, $15{y^{VI}}^2+4{y^{VI}}y^V-5{y^V}^2$, and $15{y''}^2-4y''y'''-5{y'''}^2$, $6{y^{IV}}^2-2{y^{IV}}y'''-13{y'''}^2$, $7{y^{VI}}^2+6y^{VI}y^V-10{y^V}^2$; from this and from what has already been found above I conclude that the twelve formulas that we have given for the divisors of the numbers of the form $t^2-79u^2$ may be reduced to these four:
\begin{eqnarray*}
y^2 - 79z^2, \;\;&\;\; 3y^2\pm 2yz -26z^2, \\
79z^2-y^2, \;\;&\;\; 26z^2\mp 2yz -3y^2,
\end{eqnarray*}
which must be regarded as essentially different from one another, so that they do not admit any further reduction. \\

28. According to these principles we may construct two Tables for the forms of the odd divisors of the numbers $t^2+au^2$ and $t^2-au^2$ supposing successively that $a=1$, $2$, $3$, etc. \\

Here are the Tables constructed\footnote{Literally, ``pushed.''} up to $a=31$. It would be good to continue them at least up to $100$, but we content ourselves here to put on the path those who will take charge of this work in the future. \\

We will remark, with regard to the second Table, that the ambiguous signs $\pm$ which we find, denote that the values of $p$ and $r$ which are affected may be taken equally well with the positive or negative signs; so, since $a=2$ results in $p=\pm 1$, $q=0$, $r=\pm 2$, it follows that every odd divisor of $t^2-2a^2$ will be simultaneously of the form $y^2-2z^2$ and $2z^2-y^2$, and so on. So in this case we will be free to take the positive or negative signs. \\

We must remark again that we have omitted for simplicity all the values of $a$ which are equal to squares or are divisible by squares; this is because in the column of values of $a$ we cannot find either the number $4$, nor the number $8$, etc.; indeed, it is evident that the formula $t^2+4u^2$ is included under this one, $t^2+u^2$ where $a=1$; we also see that the formula $t^2+8u^2$ is reducible to this one, $t^2+2u^2$ where $a=2$; and similarly for the others. \\

\begin{center}
\underline{\hskip 150pt}
\end{center}

\break

\begin{figure}[ht]
\centering
\textsc{Table I.}

\bigskip

Formula of the proposed numbers: $t^2+au^2$ \\

\bigskip

Formula of their odd divisors: $py^2\pm 2qyz+rz^2$, where $pr-q^2=a$. \\

\bigskip

\renewcommand{\arraystretch}{1.25}
\begin{tabular}{|c|p{70pt}|p{70pt}|p{70pt}|}
\hline
Values & \multicolumn{2}{r}{Corresponding values of} & \\
of $a$ & $p$ & $q$ & $r$ \\
\hline
1 & 1 & 0 & 1 \\
2 & 1 & 0 & 2 \\
3 & 1 & 0 & 3 \\
5 & 1, 2 & 0, 1 & 5, 3 \\
6 & 1, 2 & 0, 0 & 6, 3 \\
7 & 1 & 0 & 7 \\
10 & 1, 2 & 0, 0 & 10, 5 \\
11 & 1, 3 & 0, 1 & 11, 4 \\
13 & 1, 2 & 0, 1 & 13, 7 \\
14 & 1, 2, 3 & 0, 0, 1 & 14, 7, 5 \\
15 & 1, 3 & 0, 0 & 15, 5 \\
17 & 1, 2, 3 & 0, 1 & 17, 9, 6 \\
19 & 1, 4 & 0, 1 & 19, 5 \\
21 & 1, 3, 2, 5 & 0, 0, 1, 2 & 21, 7, 11, 5 \\
22 & 1, 2 & 0, 0 & 22, 11 \\
23 & 1, 3 & 0, 1 & 23, 8 \\
26 & 1, 2, 3, 5 & 0, 0, 1, 2 & 26, 13, 9, 6 \\
29 & 1, 3, 5 & 0, 1, 1 & 29, 10, 6 \\
30 & 1, 3, 5, 2 & 0, 0, 0, 1 & 30, 10, 6, 17 \\
31 & 1, 5 & 0, 2 & 31, 7 \\
\hline
\end{tabular}
\end{figure}

\break

\begin{figure}[ht]
\centering
\textsc{Table II.}

\bigskip

Formula of the proposed numbers: $t^2-au^2$ \\

\bigskip

Formula of their odd divisors: $py^2\pm 2qyz-rz^2$, where $pr+q^2=a$. \\

\bigskip

\renewcommand{\arraystretch}{1.25}
\begin{tabular}{|c|p{70pt}|p{60pt}|p{90pt}|}
\hline
Values & \multicolumn{2}{r}{Corresponding values of} & \\
of $a$ & $p$ & $q$ & $r$ \\
\hline
1 & 1 & 0 & 1 \\
2 & $\pm 1$ & 0 & $\pm 2$ \\
3 & 1, $-1$ & 0 & 3, $-3$ \\
5 & $\pm 1$ & 0 & $\pm 5$ \\
6 & 1, $-1$ & 0 & 6, $-6$ \\
7 & 1, $-1$ & 0 & 7, $-7$ \\
10 & $\pm 1$, $\pm 2$ & 0 & $\pm 10$, $\pm 5$ \\
11 & 1, $-1$ & 0 & 11, $-11$ \\
13 & $\pm 1$ & 0 & $\pm 13$ \\
14 & 1, $-1$ & 0 & 14, $-14$ \\
15 & 1, $-1$, 3, $-3$ & 0 & 15, $-15$, 5, $-5$ \\
17 & $\pm 1$ & 0 & $\pm 17$ \\
19 & 1, $-1$ & 0 & 19, $-19$ \\
21 & 1, $-1$ & 0 & 21, $-21$ \\
22 & 1, $-1$ & 0 & 22, $-22$ \\
23 & 1, $-1$ & 0 & 23, $-23$ \\
26 & $\pm 1$, $\pm 2$ & 0 & $\pm 26$, $\pm 13$ \\
29 & $\pm 1$ & 0 & $\pm 29$ \\
30 & 1, $-1$, 2, $-2$ & 0 & 30, $-30$, 15, $-15$ \\
31 & 1, $-1$ & 0 & 31, $-31$ \\
\hline
\end{tabular}

\bigskip\bigskip
\emph{One will find the sequel to this research in the Volume for the year 1774.}
\bigskip\bigskip
\underline{\hskip 150pt}
\end{figure}

\end{document}